\documentclass[11pt,leqno]{article}

\usepackage{amsmath,amsfonts,amscd,amssymb,theorem}

\long\def\comment#1\endcomment{}

\comment
\endcomment

\makeatletter
\begingroup
\gdef\th@dotted{\normalfont\itshape
  \def\@begintheorem##1##2{%
        \item[\hskip\labelsep \theorem@headerfont ##1\ ##2.]}%
\def\@opargbegintheorem##1##2##3{%
   \item[\hskip\labelsep \theorem@headerfont ##1\ ##2\ (##3).]}}
\endgroup
\makeatother

\theoremstyle{dotted}

\newtheorem{theorem}{Theorem}[section]
\newtheorem{lemma}[theorem]{Lemma}

\newtheorem{prop}[theorem]{Proposition}

\makeatletter
\begingroup
\gdef\th@upshape{\normalfont
  \def\@begintheorem##1##2{%
        \item[\hskip\labelsep \theorem@headerfont ##1\ ##2.]}%
\def\@opargbegintheorem##1##2##3{%
   \item[\hskip\labelsep \theorem@headerfont ##1\ ##2\ (##3).]}}
\endgroup
\makeatother

\theoremstyle{upshape}

\newtheorem{defn}[theorem]{Definition}

\makeatletter
\renewcommand{\subsection}{\@startsection{subsection}{2}{0pt}{-3ex
plus -1ex minus -0.2ex}{-2mm plus -0pt minus
-2pt}{\normalfont\bfseries}} 
\renewcommand{\subsubsection}{\@startsection{subsubsection}{3}{0pt}{-3ex
plus -1ex minus -0.2ex}{-2mm plus -0pt minus
-2pt}{\normalfont\bfseries}} 
\makeatother

\makeatletter
\@addtoreset{equation}{section}
\makeatother

\newcommand{\cntrct}                % contraction with a vector field
{\hspace{2pt}\raisebox{1pt}{\text{$\lrcorner$}}\hspace{2pt}}

\def\eqref#1{\thetag{\ref{#1}}}

\let\latexref=\ref
\def\ref#1{{\normalfont{\latexref{#1}}}}

\newcommand{\wt}{\widetilde}
\newcommand{\wh}{\widehat}

\newcommand{\idot}{{\:\raisebox{1pt}{\text{\circle*{1.5}}}}}
\newcommand{\hdot}{{\:\raisebox{3pt}{\text{\circle*{1.5}}}}}
\newcommand{\eps}{\varepsilon}
\renewcommand{\phi}{\varphi}

\newcommand{\vH}{\check{H}}

\def\dlim_#1{{\displaystyle\lim_{#1}}^\hdot}

\newcommand{\Hom}{\operatorname{Hom}}

\newcommand{\RHom}{\operatorname{RHom}}

\renewcommand{\Im}{\operatorname{Im}}

\newcommand{\id}{\operatorname{\sf id}}

\newcommand{\tr}{\operatorname{\sf tr}}

\newcommand{\A}{{\cal A}}
\newcommand{\D}{{\cal D}}
\newcommand{\DM}{{\cal D}{\cal M}}
\newcommand{\DN}{{\cal D}{\cal N}}
\newcommand{\wDM}{\widehat{\DM}}
\newcommand{\DS}{{\cal D}{\cal S}}

\newcommand{\C}{{\cal C}}
\newcommand{\B}{{\cal B}}

\newcommand{\Q}{{\cal Q}}

\newcommand{\Sets}{\operatorname{Sets}}

\newcommand{\Iso}{{\operatorname{Iso}}}
\newcommand{\Aut}{{\operatorname{Aut}}}

\newcommand{\Add}{\operatorname{\sf Add}}

\newcommand{\ppt}{{\sf pt}}

\newcommand{\bPhi}{\overline{\Phi}}

\newcommand{\M}{\operatorname{\mathcal M}}
\newcommand{\N}{\operatorname{\mathcal N}}
\newcommand{\Nn}{\operatorname{\sf N}}

\newcommand{\Infl}{\operatorname{\sf Infl}}

\newcommand{\copr}{\sqcup}

\newcommand{\Z}{{\mathbb Z}}
\newcommand{\RR}{{\mathbb R}}

\newcommand{\W}{\mathbb{W}}

\newcommand{\Fun}{\operatorname{Fun}}

\newcommand{\E}{\mathcal{E}}

\newcommand{\wM}{\operatorname{\wh{\mathcal{M}}}}
\newcommand{\wA}{\operatorname{\wh{\mathcal{A}}}}

\newcommand{\wGamma}{\operatorname{\widehat{\Gamma}}}

\newcommand{\vC}{\check{C}}

\newcommand{\DA}{\D^\alpha}

\newcommand{\hGamma}{\overline{\Gamma}}

\newcommand{\bD}{\overline{\D}}
\newcommand{\bInfl}{\overline{\Infl}}

\title{Derived Mackey functors and profunctors: an overview of results}

\author{D. Kaledin\thanks{Partially supported by RScF, grant number
    14-21-00053, by AG Laboratory SU-HSE, RF government grant,
    ag.11.G34.31.0023, and by the Dynasty Foundation award.}}

\begin{document}

\maketitle

\tableofcontents

\section*{Introduction.}

``Mackey functors'' associated to a finite group $G$ have long been
a standard tool both in group theory and in algebraic topology
(specifically, in the $G$-equivariant stable homotopy theory). The
notion was originally introduced by Dress \cite{dress} and later
clarified by several people, in particular by Lindner \cite{lind}.
The reader can find modern expositions in the topological context
e.g. in \cite{may1}, \cite{may2}, \cite{tD}, or a more algebraic
treatment in \cite{the}.

For any finite group $G$, $G$-Mackey functors form an abelian
category $\M(G)$. If one wants to develop a homological theory of
Mackey functors, it seems natural to consider the derived category
$\D(\M(G))$. However, there is an alternative to this suggested in
\cite{mackey}. By modifying the very definition of a Mackey functor,
one can construct a triangulated category $\DM(G)$ of ``derived
Mackey functors'' that contains $\M(G)$ but differs from
$\D(\M(G))$. The category $\DM(G)$ reflects better the properties of
the $G$-equivariant stable homotopy category, and it also behaves
better than $\D(\M(G))$ from the purely formal point of view. For
more details, we refer the reader to \cite{mackey}; here we only
mention two things:
\begin{itemize}
\item $\DM(G)$ has a natural semiorthogonal decomposition into
  pieces that are identified with the derived categories of
  representations of certain subquotients of $G$; the gluing
  functors of this descomposition have a natural interpretation in
  terms of a certain generalization of Tate cohomology of finite
  groups.
\item $\DM(G)$ has a lot of autoequivalences -- in particular, any
  finite-dimen\-si\-on\-al real representation $V$ of $G$ gives rise
  to a ``suspension autoequivalence'' $\Sigma^V:\DM(G) \to \DM(G)$,
  a generalization of the homological shift.
\end{itemize}
Another way to extend to the notion of a Mackey functor is to
consider more general groups $G$. In fact, if one allows groups with
non-trivial topology, then $G$ does not have to be finite -- one can
consider any compact Lie group. However, from the point of view of
algebra, a more natural modification is to allow groups that are
discrete but infinite. Formally, one can do so without any changes
-- Lindner's definition works for any group $G$, albeit the
resulting category only depends on its profinite completion
$\wh{G}$. However, once more there is an alternative suggested in
\cite{proma} under the name of a ``$G$-Mackey profunctor''. Mackey
profunctors seem to reflect better the structure of the group $G$.
In particular, in good situations, giving a $G$-Mackey profunctor
$M$ is equivalent to giving a system of Mackey functors $M_N \in
\M(G/N)$, $N \subset G$ a cofinite normal subgroup, related by some
natural compatibility isomorphisms.

One can also combine the two stories and develop the notion of a
``derived Mackey profunctor''. This has also been done in
\cite{proma}. It turns out that even if one is interested in
non-derived Mackey profunctors, looking at them in the derived
context is cleaner and gives stronger results.

\medskip

Unfortunately, both \cite{mackey} and \cite{proma} are rather long
and technical papers. The goal of the present paper is to give a
brief overview of them. We do not give any proofs whatsoever, and we
keep the technicalities to an absolute minimum. Instead, we try to
present the conceptual ideas behind the constructions, and to
illustrate the theory by some key examples.

\medskip

The paper is organized as follows. Section 1 is a recollection of
the standard theory of Mackey functors. Section 2 is concerned with
derived Mackey functors for a finite group $G$. Generally, we follow
\cite{mackey}, but we use improvements suggested in \cite{proma} to
clean up some statements and simplify the proofs. In particular, we
describe and use a relation between Mackey functors and finite
pointed $G$-sets. We illustrate the theory in the simplest possible
example $G=\Z/p\Z$, the finite cyclic group of prime order. In
Section 3, we turn to the theory of Mackey profunctors developed in
\cite{proma}. The main example here is $G=\Z$, the infinite cyclic
group.

\subsection*{Acknowledgements.} I am grateful to the organizers of
the International School on TQFT, Langlands and Mirror Symmetry in
Playa del Carmen, Mexico, in March 2014 for inviting me and giving
me the opportunity to present my results. It is a pleasure to submit
this overview to the proceedings.

\section*{Notation.}

For any category $\E$, we denote by $\E^o$ the opposite
category. For every integer $n \geq 0$, we denote by $[n]$ the set
of integers $i$, $0 \leq i \leq n$, and as usual, we denote by
$\Delta$ the category of the sets $[n]$ and order-preserving maps
$g$ between them. Simplicial objects in a category $\E$ are functors
from $\Delta^o$ to $\E$.

For any small category $I$ and ring $R$, we denote by $\Fun(I,R)$
the abelian category of functors from $I$ to the category of
$R$-modules, and we denote by $\D(I,R)$ the derived category of the
abelian category $\Fun(I,R)$. For any functor $\gamma:I' \to I$
between small categories, we have the natural pullback functor
$\gamma^*:\Fun(I,R) \to \Fun(I',R)$, and it has a left and a right
adjoint $\gamma_!,\gamma_*:\Fun(I',R) \to \Fun(I,R)$, the left and
right Kan extension functors. The derived functors $L^\hdot
\gamma_!,R^\hdot\gamma_*:\D(I',R) \to \D(I,R)$ are left and
right-adjoint to the pullback functor $\gamma^*:\D(I,R) \to
\D(I',R)$.

\section{Recollection on Mackey functors.}\label{mack.sec}

We start with a brief recapitulation of the usual theory of Mackey
functors (we follow the exposition in \cite[Section 2]{proma}).

Assume given a category $\C$ that has fibered products. The {\em
  category $Q\C$} is the category wwith same objects as $\C$, and
with morphisms from $c_1$ to $c_2$ given by isomorphism classes of
diagrams
\begin{equation}\label{domik}
\begin{CD}
c_1 @<{l}<< c @>{r}>> c_2
\end{CD}
\end{equation}
in $\C$. The compositions are obtained by taking
pullbacks. Note that we have a natural embedding
\begin{equation}\label{e.eq}
e:\C^o \to Q\C
\end{equation}
that is identical on objects, and sends a map $f$ to the diagram
\eqref{domik} with $l=f$ and $r=\id$. The category $Q\C$ is
obviously self-dual: $Q\C^o$ is identified with $Q\C$ by the functor
that flips $r$ and $l$ in \eqref{domik}. By duality, the embedding
\eqref{e.eq} induces an embedding $e^o:\C \to Q\C$.

Now fix a group $G$ and a ring $R$, and take as $\C$ the category
$\Gamma_G$ of finite $G$-sets, that is, finite sets equipped with an
action of the group $G$. The category $\Gamma_G$ has fibered
products, so we can define the category $Q\Gamma_G$ and the
embedding \eqref{e.eq}.

\begin{defn}\label{add.def}
An object $E \in \Fun(\Gamma_G^o,R)$ is {\em additive} if for any
$S_1,S_2 \in \Gamma_G^o$, the natural map
\begin{equation}\label{add.eq}
E(S_1 \copr S_2) \to E(S_1) \oplus E(S_2)
\end{equation}
is an isomorphism. An {\em $R$-valued $G$-Mackey functor} $E$ is an
object $E \in \Fun(Q\Gamma_G,R)$ whose restriction $e^*E \in
\Fun(\Gamma_G^o,R)$ is additive. The full subcategory in
$\Fun(Q\Gamma_G,R)$ formed by Mackey functors is denoted by
$\M(G,R)$.
\end{defn}

Denote by $\Fun_{add}(\Gamma_G^o,R) \subset \Fun(\Gamma_G^o,R)$ the
full subcategory spanned by additive objects, and let $O_G \subset
\Gamma_G$ be the subcategory of $G$-orbits -- that is, $G$-sets of
the form $[G/H]$, $H \subset G$ a subgroup. Since every finite
$G$-set decomposes into a disjoint union of $G$-orbits,
\eqref{add.eq} insures that restriction to $O_G^o \subset
\Gamma_G^o$ provides an equivalence
$$
\Fun_{add}(\Gamma_G^o,R) \cong \Fun(O_G^o,R).
$$
In other words, an additive object $E \in \Fun(\Gamma_G^o,R)$ is
uniquely determined by specifyings its values $E^H = E([G/H])$ at all
$G$-orbits, and maps
$$
f^* = E(f):E^H \to E^{H'},
$$
one for each map $f:[G/H'] \to [G/H]$ of $G$-orbits, compatible with
compositions for composable pairs of maps $f$, $f'$. To extend $E$
to a Mackey functor, one needs to also specify a map
\begin{equation}\label{trans}
f_* = E(e^o(f)):E(S') \to E(S)
\end{equation}
for any map $f:S' \to S$ of finite $G$-sets. These are sometimes
called {\em transfer maps} (this is how they appear in equivariant
stable homotopy theory). By additivity, it suffices to specify the
transfer maps for maps $G$-orbits. They should be compatible with
compositions, and for any two maps $f':S' \to S$, $f'':S'' \to S$,
there is a certain compatibility condition between $f''_*$ and
$f^{'*}$ encoded in the definition of the category
$Q\Gamma_G$. Explicitly, if $S=[G/H]$, $S' = [G/H']$, $S'' =
[G/H'']$, and $f'$, $f''$ are induced by embeddings $H',H'' \subset
H$, the fibered product $S' \times_S S''$ decomposes into a disjoint
union
\begin{equation}\label{orb.spl}
S' \times_S S'' = \coprod_{s \in S}[G/H_s]
\end{equation}
of $G$-orbits indexed by the finite set $S = H' \setminus\! H / H''$,
and we must have
\begin{equation}\label{double}
f^{'*} \circ f''_* = \sum_{s \in S}\wt{f}''_{s*} \circ
\wt{f}^{'*}_s,
\end{equation}
where $f'_s:[G/H_s] \to [G/H'']$, $f''_s:[G/H_s] \to [G/H']$ are
projections of the component $G/H_s$ of the decomposition
\eqref{orb.spl}. This is known as the {\em double coset formula}.

For any subgroup $H \subset G$, we have a natural restriction
functor $\rho^H:\Gamma_G \to \Gamma_H$ sending a $G$-set $S$ to the
same set $S$ considered as an $H$-set. The functor $\rho^H$ has a
left-adjoint $\gamma^H:\Gamma_H \to \Gamma_G$. Both $\rho^H$ and
$\gamma^H$ preserve fibered products, thus induce functors
\begin{equation}\label{q.ga}
Q(\gamma^H):Q\Gamma_H \to Q\Gamma_G, \qquad Q(\rho^H):Q\Gamma_G \to
Q\Gamma_H.
\end{equation}
One checks that the adjunction between $\rho^H$ and $\gamma^H$
induces an adjunction between $Q(\gamma^H)$ and $Q(\rho^H)$, and
that the functor
\begin{equation}\label{psi.h}
\Psi^H = Q(\rho^H)_! \cong Q(\gamma^H)^*:\Fun(Q\Gamma_G,R) \to
\Fun(Q\Gamma_H,R)
\end{equation}
preserves additivity, thus sends Mackey functors to Mackey
functors. The functor $\Psi^H$ is known as the {\em categorical
  fixed points functor}.

On the other hand, for any subgroup $H \subset G$ with normalizer
$N_H \subset G$, the quotient $W = N_H/H = \Aut_G([G/H])$ acts on
the fixed points subset $S^H$ of any $G$-set $S$. Then sending $S$
to $S^H$ gives a functor $\phi^H:\Gamma_G \to \Gamma_W$ preserving
fibered products, and we can consider the adjoint pair of functors
\begin{equation}\label{infl.h}
\begin{aligned}
\Phi^H &= Q(\phi^H)_!:\Fun(Q\Gamma_G,R) \to \Fun(Q\Gamma_W,R),\\
\Infl^H &= Q(\phi^H)^*:\Fun(Q\Gamma_W,R) \to \Fun(Q\Gamma_G,R).
\end{aligned}
\end{equation}
Both preserve additivity, thus induce adjoint functors between
$\M(G,R)$ and $\M(W,R)$. These are the {\em geometric fixed points
  functor} and the {\em inflation functor}.

The embedding $\M(G,R) \subset \Fun(Q\Gamma_G,R)$ admits a
left-adjoint {\em additivization functor} $\Add$. This allows to
define tensor products of Mackey functors. Namely, the cartesian
product functor
\begin{equation}\label{m.G}
m:\Gamma_G \times \Gamma_G \to \Gamma_G
\end{equation}
preserves fiber products, thus induces a functor
$$
Q(m):Q\Gamma_G \times Q\Gamma_G \to Q\Gamma_G.
$$
Then for any two algebras $R_1$, $R_2$ over a commutative ring $k$
and Mackey functors $E_1 \in \M(G,R_1)$, $E_2 \in \M(G,R_2)$, their
{\em tensor product} $E_1 \circ E_2$ is given by
\begin{equation}\label{ma.produ}
E_1 \circ E_1 = \Add(Q(m)_!(E_1 \boxtimes_k E_2)) \in \M(G,R_1 \otimes_k
R_2).
\end{equation}
The tensor product \eqref{ma.produ} is unital, with the unit given
by the so-called {\em Burnside Mackey functor} $\A \in \M(G,\Z)$.
Explicitly, for any subgroup $H \subset G$, $\A^H$ is canonically
identified with the {\em Burnside ring} $A^H$ of the group $H$ given
by
\begin{equation}\label{A.H.eq}
A^H = \Z[\Iso(\Gamma_H)]/\{[S \copr S'] - [S] - [S']\},
\end{equation}
that is, the free abelian group generated by isomorphism classes
$[S]$ of finite $H$-sets $S$, modulo the relations $[S \copr S'] =
[S] + [S']$ (this is a commutative ring, with the product induced by
the cartesian product of $H$-sets). Note that since any $G$-set
uniquely decomposes into a disjoint union of $G$-orbits, we have
\begin{equation}\label{a.g}
A^G \cong \bigoplus_{H \subset G}\Z,
\end{equation}
where the sum is over all conjugacy classes of subgroups $H \subset
G$. Both the categorical fixed points functor $\Psi^H$ and the
geometric fixed points functor $\Phi^H$ commute with the product
\eqref{ma.produ}.  If $R$ is a commutative ring, then the product
\eqref{ma.produ} with $R_1 = R_2 = k = R$ turns $\M(G,R)$ into a
unital symmetric tensor category. An {\em $R$-valued $G$-Green
  functor} is an algebra object in the category $\M(G,R)$.

For any subgroup $H \subset G$, the functor $\gamma^H$ and the
categorical fixed points functor $\Psi^H$ can be computed rather
explicitly. The geometric fixed points functor $\Phi^H$ is harder to
describe. However, there is the following result.

\begin{lemma}[{{\cite[Lemma 2.4]{proma}}}]\label{infl.le}
Assume given a normal subgroup $N \subset G$, with the quotient $W =
G/N$. Then the inflation functor $\Infl^N$ is fully
faithful. Moreover, for any $E \in \M(G,R)$, the adjunction map
\begin{equation}\label{infl.adj}
M \to \Infl^N\Phi^NM
\end{equation}
is surjective, and for any $S \in \Gamma_W$, we have a short
exact sequence
\begin{equation}\label{phi.exp}
\begin{CD}
\displaystyle\bigoplus_{f:S' \to S}M(S') @>{\sum f_*}>> M(S) @>>>
\Phi^N(M)(S) @>>> 0,
\end{CD}
\end{equation}
where $f_*$ is as in \eqref{trans}, and the sum is over all maps
$f:S' \to S$ in $\Gamma_G$ such that $S'$ has no elements fixed
under $N \subset G$.
\end{lemma}

Note that by additivity of $E \in \Fun(Q\Gamma_G,R)$, the image of
the map $\sum f_*$ in \eqref{phi.exp} is not only the sum of the
images of individual maps $f_*$, it is actually the union of these
images.

\section{Derived Mackey functors.}\label{dm.sec}

We now turn to the derived version of the story of
Section~\ref{mack.sec}. To see why the derived category
$\D(\M(G,R))$ might be a wrong object to work with, it suffices to
recall the definition of the category $Q\C$: we only consider the
isomorphism classes of diagrams \eqref{domik}, and completely ignore
the fact that those diagrams can have non-trivial automorphism
groups with non-trivial homology. On the formal level, the problem
appears for example in Lemma~\ref{infl.le}: if one considers the
derived categories $\D(\M(G,R))$, $\D(\M(W,R))$, then the functor
$$
\Infl^N:\D(\M(W,R)) \to \D(\M(G,R))
$$
induced by the inflation functor \eqref{infl.h} is not longer fully
faithful. To obtain the correct category of derived Mackey functors,
one has to take account of the automorphisms of diagrams
\eqref{domik} and treat $Q\Gamma_G$ as what it naturally is -- a
$2$-category rather than simply a category. But then, one has to
make sense of the derived category $\D(Q\Gamma_G,R)$.

\subsection{Quotient constructions.}\label{sc.subs}

In fact, \cite{mackey} introduces not one but two equivalent ways to
do it. Firstly, one can observe that any small category $\C$ defines
a small additive category $\B^\C$ with the same objects, and
morphism groups given by
$$
\B^\C(c,c') = \Z[\C(c,c')],
$$
where $\C(c,c')$ is the set of morphisms from $c$ to $c'$ in the
category $\C$. Giving a functor from $\C$ to $R$-modules is
equivalent to giving an additive functor from $\B^\C$ to
$R$-modules. If $\C$ is not a category but a $2$-category, then
$\C(c,c')$ are not sets but small categories. What one can do is
consider the geometric realizations $|\C(c,c')|$ of the nerves of
the categories $\C(c,c')$, and let
$$
\B^\C_\idot = C_\idot(|\C(c,c')|,\Z)
$$
be the singular chain complexes of these geometric realizations. If
$\C$ is a strict $2$-category -- that is, compositions are
associative on the nose -- then geometric realization is
sufficiently functorial so that $\B^\C_\idot(-,-)$ define a DG
category $\B^\C_\idot$ with the same objects as $\C$. In the general
case, it is still possible to define $\B^\C_\idot$ as an
$A_\infty$-category. Then one can define $\D(\C,R)$ as the dervied
category of $A_\infty$-functors from $\B^\C_\idot$ to complexes of
$R$-modules. This is done in \cite[Subsection 1.6]{mackey}.

The second approach is that of \cite[Section 4]{mackey}. It uses the
nerves more directly. Namely, recall that the nerve $N(\C):\Delta^o
\to \Sets$ of a small category $\C$ is a simplicial set whose value
at $[n] \in \Delta^o$ is the set of all diagrams
\begin{equation}\label{n-ex}
\begin{CD}
c_0 @>{f_0}>> \dots @>{f_{n-1}}>> c_n
\end{CD}
\end{equation}
in $\C$. For any map $g:[n] \to [n']$ in $\Delta$, a diagram
$c_\idot \in N(\C)([n'])$ induces a diagram $g^*c_\idot \in
N(\C)([n])$ such that $g^*c_i = c_{g(i)}$. One can turn $N(\C)$ into
another small category $\N(\C)$ by applying the Grothendieck
construction of \cite{SGA}: the objects of $\N(\C)$ are pairs
$\langle [n],c_\idot \rangle$ of an object $[n] \in \Delta$ and a
diagram \eqref{n-ex}, and morphisms from $\langle [n],c_\idot
\rangle$ to $\langle [n'],c'_\idot \rangle$ are morphisms $g:[n] \to
[n']$ in $\Delta$ such that $c_\idot = g^*c'_\idot$. Sending a
diagram \eqref{n-ex} to $c_n$ gives a functor $q:\N(\C) \to
\C$, and the corresponding pullback functor
$$
q^*:\D(\C,R) \to \D(\N(\C),R)
$$
is a fully faithful embedding. To characterize its image, say that a
morphism in $\N(\C)$ defined by $g:[n] \to [n']$ is {\em special} if
$g(n) = n'$, and say that $E \in \D(\N(\C),R)$ is {\em special} if
the map $E(f)$ is invertible for any special $f$. Then $\D(\C,R)
\subset \D(\N(\C),R)$ is exactly the full subcategory spanned by
special objects.

Now, if one starts with a $2$-category $\C$, then diagrams
\eqref{n-ex} form a category rather than a set. However, one can
still define the category $\N(\C)$: objects are pairs $\langle
[n],c_\idot \rangle$, morphisms from $\langle [n],c_\idot \rangle$
to $\langle [n'],c'_\idot \rangle$ are pairs of a morphism $g:[n] \to
[n']$ and a morphism $c_\idot \to g^*c'_\idot$. Then one can keep the
same notion of a special map and special object, and define
$\D(\C,R)$ as the full subcategory in $\D(\N(\C),R)$ spanned by
special objects.

As it turns out, the second approach is much better for practical
applications: one does not need to keep track of all the higher
$A_\infty$-operations, since they are all packed into the structure
of the category $\N(\C)$.

\medskip

For applications to Mackey functors, then, one would start with a
category $\C$ that has fibered products, define a category
$Q\C$ with the same objects and morphisms given by isomorphisms
classes of diagrams
\begin{equation}\label{domik.bis}
\begin{CD}
c @<<< \wt{c} @>>> c'
\end{CD}
\end{equation}
in $\C$, and refine it to a $2$-category $\Q\C$ by letting
$\Q\C(c,c')$ be the groupoids of such diagrams and invertible maps
between them. Then one would consider the nerve $\N(\Q\C)$. This is
more-or-less what is done in \cite[Section 4]{mackey} and
\cite[Section 4]{proma}, but with a small modification. It turns out
that instead of $\N(\Q\C)$, one can use a smaller category
$S\C$. Its objects are pairs $\langle [n],c_\idot \rangle$ of $[n]
\in \Delta$ and a diagram \eqref{n-ex} in $\C$, and morphisms from
$\langle [n],c_\idot \rangle$ to $\langle [n'],c'_\idot \rangle$ are
pairs $\langle g,\alpha \rangle$ of a morphism $g:[n] \to [n']$ and
a collection of morphisms $\alpha_i:c'_{g(i)} \to c_i$, $0 \leq i
\leq n$, such that for any $i$ and $j$, $1 \leq i \leq j \leq n'$,
the diagram
$$
\begin{CD}
c'_{g(i)} @>>> c'_{g(j)}\\
@V{\alpha_i}VV @VV{\alpha_j}V\\
c_i @>>> c_j
\end{CD}
$$
is a fibered product square in the category $\C$. A map $\langle
g,\alpha_\idot \rangle$ in $S\C$ is special if $g(n) = n'$ and
$\alpha_n$ is invertible, and an object $E \in \D(S\C,R)$ is special
if $E(f)$ is invertible for any special map $f$. One denotes by
\begin{equation}\label{ds.eq}
\DS(\C,R) \subset \D(S\C,R)
\end{equation}
the full subcategory spanned by special objects.

One advantage of the category $S\C$ is that one can define an analog
of the embedding \eqref{e.eq} without replacing the category $\C^o$
with its nerve: by definition, sending $c \in \C$ to the pair
$\langle [0],c \rangle$ gives a natural functor $e:\C^o \to
S\C$. Thus every object $E \in \D(S\C,R)$ defines an object
\begin{equation}\label{bp.eq}
\overline{E} = e^*E \in \D(\C^o,R)
\end{equation}
that we call the {\em base part} of $E$. Sending a diagram
\eqref{n-ex} to $c_n$ gives a natural functor $q:S\C \to Q\C$, and
the composition $q \circ e$ is the natural embedding \eqref{e.eq}.
For every special map $f$ in $S\C$, $q(f)$ is invertible in $Q\C$,
so that we have a natural functor
\begin{equation}\label{q.st}
q^*:\D(Q\C,R) \to \DS(\C,R).
\end{equation}
This functor is in general not an equivalence, and we take
$\DS(\C,R)$ as the correct definition of $\D(\Q\C,R)$, a refinement
of $\D(Q\C,R)$.

Although $q^*$ is not an equivalence, we note that the standard
$t$-structure on $\D(S\C,R)$ induces a $t$-structure on $\DS(\C,R)$,
$q^*$ is $t$-exact, and it identifies the heart of the $t$-structure
on $\DS(\C,R)$ with $\Fun(Q\C,R)$. For any $E \in \Fun(Q\C,R)$, the
base part of $q^*E$ is $e^*E \in \Fun(\C^o,R)$; extending this base
part to a special object in $\D(S\C,R)$ is equivalent to providing
the transfer maps \eqref{trans}.

The category $\DS(\C,R)$ of \eqref{ds.eq} is functorial in $\C$ in
the following sense. Say that a functor $F:\C \to \C'$ between two
small categories with fibered product is a {\em morphism} if it
preserves fibered products. Then any morphism $F$ induced a functor
$S(F):S\C \to S\C'$ that commutes with projections to $\Delta$ and
sends special maps to special maps. Therefore we have a pullback
functor
$$
S(F)^*:\DS(\C',R) \to \DS(\C,R).
$$
It has been shown in \cite[Corollary 4.15]{proma} that $S(F)^*$ has
a left-adjoint functor
$$
S(F)_!:\DS(\C,R) \to \DS(\C,R).
$$

\subsection{Definitions and properties.}

Now return to the situation of Section~\ref{mack.sec}: fix a finite
group $G$, and let $\Gamma_G$ be the category of finite
$G$-sets. This category has fiberd products, so that for any ring
$R$, we can consider the derived category $\DS(\Gamma_G,R)$ of
\eqref{ds.eq}.

\begin{defn}
An object $E \in \D(\Gamma_G^o,R)$ is {\em additive} if the natural
map \eqref{add.eq} is an isomorphism for any $S_1,S_2 \in
\Gamma_G$. A {\em derived $R$-valued $G$-Mackey functor} is an
object $E \in \DS(\Gamma_G,R)$ with additive base part
\eqref{bp.eq}.
\end{defn}

We denote the full subcategory of derived $R$-valued $G$-Mackey
functors by $\DM(G,R) \subset \DS(\Gamma_G,R)$. The standard
$t$-structure on $\DS(\Gamma_G,R)$ induces a $t$-structure on
$\DM(G,R)$ whose heart is identified with $\M(G,R)$ by the functor
$q^*$ of \eqref{q.st}. It has been proved in \cite[Lemma 8.3]{proma}
that the embedding $\DM(G,R) \subset \DS(\Gamma_G,R)$ admits a
left-adjoint {\em additivization functor}
\begin{equation}\label{Add.der}
\Add:\DS(\Gamma_G,R) \to \DM(G,R).
\end{equation}
If $G = \{e\}$ is the trivial group consisting of its unity element
$e$, then $\DM(G,R) \cong \D(R)$ is the derived category of the
category of $R$-modules (note that this is a non-trivial statement
that requires a proof, such as the one found in \cite[Subsection
  3.2]{mackey}). As in \eqref{q.ga}, for any subgroup $H \subset G$,
we have the adjoint pair of functors $\rho^H$, $\gamma^H$; both are
morphisms and induce an adjoint pair of functors
$$
S(\gamma^H):S\Gamma_H \to S\Gamma_G, \qquad S(\rho^H):S\Gamma_G \to
S\Gamma_H.
$$
The functor
$$
S(\gamma^H)^* \cong S(\rho^H)_!:\DS(\Gamma_G,R) \to \DS(\Gamma_H,R)
$$
sends additive objects to additive objects, thus induces a functor
\begin{equation}\label{psi.h.der}
\Psi^H:\DM(G,R) \to \DM(H,R).
\end{equation}
This is the derived counterpart of the categorical fixed points
functor \eqref{psi.h}. On the other hand, for any $H \subset G$ with
normalizer $N_H$ and quotient $W = N_H/H$, the fixed points functor
$\phi^H:\Gamma_G \to \Gamma_W$ is also a morphism, and induces an
adjoint pair of functors
\begin{equation}\label{infl.h.der}
\begin{aligned}
\Phi^H &= S(\phi^H)_!:\DS(\Gamma_G,R) \to \DS(\Gamma_W,R),\\
\Infl^H &= S(\phi^H)^*:\DS(\Gamma_W,R) \to \DS(\Gamma_G,R),
\end{aligned}
\end{equation}
a derived counterpart the geometric fixed points functor and the
inflation functor of \eqref{infl.h}. As in the non-derived case,
both $\Phi^H$ and $\Infl^H$ preserve additivity, thus induce an
adjoint pair of functors between $\DM(G,R)$ and $\DM(H,R)$ (for
$\Phi^H$, this is \cite[Lemma 6.14]{proma}). As in
Lemma~\ref{infl.le}, for any normal subgroup $N \subset G$ with
quotient $W=G/N$, the inflation functor $\Infl^N:\DM(W,R) \to
\DM(G,R)$ is fully faithful (this is \cite[Lemma 6.12]{proma}).

To define products of derived Mackey functors, one considers the
cartesian product functor \eqref{m.G}. It is a morphism, but
$S(m)_!$ does not automatically preserve additivity. So, for any two
algebra $R_1$, $R_2$ over a commutative ring $k$ and any $E_1 \in
\DM(G,R_1)$, $E_2 \in \DM(G,R_2)$, we let
\begin{equation}\label{ma.produ.der}
E_1 \circ E_2 = \Add(S(m)_!(E_1 \boxtimes_k E_2)) \in \DM(G,R_1 \otimes_k
R_2),
\end{equation}
where $\Add$ is the additivization functor \eqref{Add.der}. This
product has exactly the same properties as the underived product
\eqref{ma.produ}. In particular, the fixed points functors $\Psi^H$,
$\Phi^H$ are tensor functors. Moreover, we have a derived version
$\A_\idot \in \DM(G,\Z)$ of the Burnside Mackey functor that serves
as the unit object for the product, and a derived Burnside ring
$A^G_\idot = \A_\idot([G/G])$. Explicitly, as shown in
\cite[Subsection 8.3]{proma}, we have
$$
A^G_\idot = \bigoplus_{H \subset G}C_\idot(W_H,\Z),
$$
where the sum is over all the conjugacy classes of subgroups $H
\subset G$, and the corresponding summand is the homology complex of
the group $W_H = \Aut_G([G/H])$ with coefficients in $\Z$. In
homological degree $0$, this recovers the isomorphism \eqref{a.g}.

\subsection{Pointed $G$-sets.}

To obtain more information about derived Mackey functors, it is
convenient to use the following observation. Say that a diagram
\eqref{domik} is {\em restricted} if the map $l$ is injective, and
note that the pullback of an injective map is an injective
map. Therefore we have a subcategory $Q_I\Gamma_G \subset Q\Gamma_G$
whose maps are isomorphism classes of restricted diagrams, and a
subcategory $S_I\Gamma_G \subset \Gamma_G$ whose maps are pairs
$\langle g,\alpha_\idot \rangle$ with injective $\alpha_n$. Say that
a map in $S_I\Gamma_G$ is special if it is special in $S\Gamma_G$,
say that an object $E \in \D(S_I\Gamma_G,R)$ is special if it
inverts all special maps, and let $\DS_I(\Gamma_G,R) \subset
\D(S_I\Gamma_G,R)$ be the full subcategory spanned by special
objects. Then the projection $q:S\Gamma_G \to Q\Gamma_G$ restricts
to a projection $q:S_I\Gamma_G \to Q_I\Gamma_G$, and we have a
natural functor
\begin{equation}\label{q.st.bis}
q^*:\D(Q_I\Gamma_G,R) \to \DS_I(\Gamma_G,R).
\end{equation}
However, unlike the functor \eqref{q.st}, this functor is an
equivalence of categories (this is \cite[Corollary 4.17]{proma} ---
roughly speaking, the reason this holds is that restricted diagrams
\eqref{domik} have no non-trivial automorphisms).

The category $Q_I\Gamma_G$ is naturally identified with the category
$\Gamma_{G+}$ of finite pointed $G$-sets (that is, finite $G$-sets
$S$ with distinguished $G$-invariant element $o \in S$). The
equivalence sends $S$ to $S \setminus \{o\}$, and a map $f:S \to S'$
goes to the diagram
$$
\begin{CD}
S \setminus \{o\} @<<< S \setminus f^{-1}(\{o'\}) @>{f}>> S'
\setminus \{o'\}.
\end{CD}
$$
The equivalence \eqref{q.st.bis} together with the restriction with
respect to the embedding $S_I\Gamma_G \to S\Gamma_G$ then provide a
natural functor
\begin{equation}\label{res}
r:\DS(\Gamma_G,R) \to \D(\Gamma_{G+},R),
\end{equation}
and we introduce the following definition.

\begin{defn}[{{\cite[Definition 6.9]{proma}}}]
For any derived Mackey functor $E \in \DM(G,R)$ and any simplicial
finite pointed $G$-set $X:\Delta^o \to \Gamma_{G+}$, the {\em
  homology complex} of $X$ with coefficients in $E$ is given by
$$
C_\idot(X,E) = C_\idot(\Delta^o,X^*r(E)),
$$
where $r$ is the restriction functor \eqref{res}, and
$C_\idot(\Delta^o,-)$ is the homology complex of the small category
$\Delta^o$.
\end{defn}

Moreover, in \cite[Subsection 7.4]{proma}, this definition is
extended in the following way: for any $E \in \DM(G,R)$ and $X \in
\Delta^o\Gamma_{G+}$, one defines a product $X \wedge E \in
\DM(G,R)$. This product is functorial in $X$ and $E$. For any $S \in
\Gamma_G$, one has a natural identification
\begin{equation}\label{wedge.circ}
(X \wedge E)(S) \cong C_\idot(X \wedge S_+,E),
\end{equation}
where $S_+ \in \Gamma_{G+}$ is obtained by adding a distinguished
element $o$ to the set $S$, and $X \wedge S_+$ stands for pointwise
smash-product of pointed sets. With this product, one can prove the
following analog of \eqref{phi.exp}.

\begin{defn}
A simplicial finite pointed $G$-set $X \in \Delta^o\Gamma_{G+}$ is
{\em adapted} to a normal subgroup $N \subset G$ if
\begin{enumerate}
\item $X^N = [1]_+$ is the constant pointed simplicial set with two
  elements, one distinguished, one not, and
\item for any subgroup $H \subset G$ not containing $N$, the reduced
  chain homology complex $\overline{C}_\idot(X^H,\Z)$ is acyclic.
\end{enumerate}
\end{defn}

\begin{lemma}[{{\cite[Lemma 7.14]{proma}}}]\label{infl.der.le}
Assume given a finite pointed simplicial $G$-set $X_N \in
\Delta^o\Gamma_{G+}$ adapted to a normal subgroup $N \subset
G$. Then for any $E \in \DM(G,R)$, we have a natural isomorphism
$$
X \wedge E \cong \Infl^N(\Phi^N(E)).
$$
\end{lemma}

We note that together with \eqref{wedge.circ},
Lemma~\ref{infl.der.le} provides a canonical identification
\begin{equation}\label{phi.exp.der}
\Phi^N(E)(S) \cong C_\idot(X_N \wedge S_+,E)
\end{equation}
for any $S \in \Gamma_W$, $W = G/N$. It is not difficult to show
that for any $N \subset G$, an adapted set $X_N \in
\Delta^o\Gamma_{G+}$ does exist. For example, one can take the union
$S_N$ of all $G$-orbits $[G/H]$, $H \subset G$ not containing $N$,
consider the pointed simplicial set $E(S_N)$ given by
\begin{equation}\label{E.N}
E(S_N)([n]) = S_N^{n+1},
\end{equation}
and let $X_N$ be the cone of the natural map $E(X_N) \to [1]_+$. If
for this choice of $X_N$, one computes the homology of the
simplicial category $\Delta^o$ by the standard complex, then
\eqref{phi.exp.der} gives \eqref{phi.exp} in homological degree $0$.

\begin{defn}
A simplicial finite pointed $G$-set $X \in \Delta^o\Gamma_{G+}$ is a
{\em homological sphere} if for any subgroup $H \subset G$, we have
$\overline{C}_\idot(X^H,\Z) \cong \Z[d_H]$ for some integer $d_H
\geq 0$.
\end{defn}

\begin{lemma}[{{\cite[Proposition 7.18]{proma}}}]\label{sph.le}
For any homological sphere $X$, the functor
$$
E \mapsto X \wedge E
$$
is an autoequivalence of the category $\DM(G,R)$.
\end{lemma}

Constructing homological spheres is as easy as constructing adapted
sets. For example, for any finite-dimensional real representation
$V$ of the group $G$, we can take a $G$-invariant triangulation of
its one-point compactification $S_V$ and obtain a homological sphere
(the dimensions $d_H$ are given by $d_H = \dim_{\RR} V^H$). Thus
Lemma~\ref{sph.le} produces many autoequivalences of the category
$\DM(G,R)$.

\subsection{An example.}

To illustrate the theory of derived Mackey functors on a concrete
example, let us consider the case $G = \Z/p\Z$, the cyclic group of
a prime order. Then there are exactly two subgroups in $G$, the
trivial subgroup $\{e\} \subset G$ and $G$ itself, so that up to an
isomorphism, $O_G$ has two objects: the free orbit $[G/\{e\}]$ and
the trivial orbit $[G/G]$. Thus a $G$-Mackey functor $E \in \M(G,R)$
is defined by two $R$-modules, $E^0 = E([G/\{e\}])$ and $E^1 =
E([G/G])$. The module $E^0$ carries the action of the group
$G$. Equivalently, if we denote by $\sigma \in \Z/p\Z$ the
generator, then we have an automorphism $\sigma:E^0 \to E^0$ such
that $\sigma^p=\id$. Moreover, there one non-trivial map
$f:[G/\{e\}] \to [G/G]$. Since $f \circ \sigma = f$, the
corresponding maps $f_*$ and $f^*$ induce natural maps
\begin{equation}\label{vf.ma}
V = f_*:(E^0)_\sigma \to E_1, \qquad F = f^*:E^1 = (E^0)^\sigma,
\end{equation}
where $(E^0)_\sigma$, $(E^0)^\sigma$ stands for coinvariants and
invariants with respect to the automorphism $\sigma$. The double
coset formula \eqref{double} then reads as
\begin{equation}\label{fv.co}
F \circ V = \id + \sigma + \dots + \sigma^{p-1}.
\end{equation}
These are the only conditions: the category $\M(G,R)$ is equivalent
to the category of pairs $\langle E^0,E^1 \rangle$ of two
$R$-modules equipped with an automorphism $\sigma:E^0 \to E^0$ of
order $p$ and two maps \eqref{vf.ma} satisfying \eqref{fv.co}.

To obtain a similar description of the category $\DM(G,R)$, choose a
resolution $P_\idot$ of the trivial $\Z[G]$-module $\Z$ by finitely
generated projective $\Z[G]$-modules, and denote
$$
C_\idot(G,-) = P_\idot \otimes_{\Z[G]} -, \qquad C^\hdot(G,-) =
\Hom_{\Z[G]}(P_\idot,-)
$$
the homology and cohomology complexes of the group $G$ computed
using $P_\idot$ (for example, one can take the standard periodic
resolution, and this would give the standard periodic
complexes). For any $R[G]$-module $E$, we have $H_0(G,E) = E_\sigma$
and $H^0(G,E) = E^\sigma$, and the functorial trace map $\id +
\sigma + \dots + \sigma^{p-1}:E_\sigma \to E^\sigma$ induces a
functorial trace map
\begin{equation}\label{tr.eq}
\tr:C_\idot(G,-) \to C^\hdot(G,-).
\end{equation}
Then $\DM(G,R)$ is obtained by inverting quasiisomorphisms in the
category of pairs $\langle E^0_\idot,E^1_\idot \rangle$ of complexes
of $R$-modules equipped with an automorphism $\sigma:E^0_\idot \to
E^0_\idot$ of order $p$ and two maps
\begin{equation}\label{f.v.dia}
\begin{CD}
C_\idot(G,E^0_\idot) @>{V}>> E^1 @>{F}>> C^\hdot(G,E^0_\idot)
\end{CD}
\end{equation}
whose composition $F \circ V$ coincides with the trace map
\eqref{tr.eq}, $F \circ V = \tr$.

This description shows clearly why the category $\DM(G,R)$ is
different from the derived category $\D(\M(G,R))$. Indeed, objects
in $\D(\M(G,R))$ are also represented by pairs $\langle
E^0_\idot,E^1_\idot \rangle$ and maps $\sigma$, $F$, $V$, but the
homology and cohomology complexes in \eqref{f.v.dia} are replaces
with coinvariants $(E^0_\idot)_\sigma$ and invariants
$(E^0_\idot)^\sigma$. To make these quasiisomorphic to the whole
homology complexes, one needs to choose a representative $E^0_\idot$
that is both $h$-projective and $h$-injective as a complex of
$R[G]$-modules. For general complex of $R[G]$-modules, such a
representative does not exists.

For any object $E \in \DM(G,R)$ represented by a pair $\langle
E^0_\idot,E^1_\idot \rangle$ as above, the geometric fixed points
$\Phi^G(E) \in \D(R)$ can be computed by
Lemma~\ref{infl.der.le}. The result is, $\Phi^G(E)$ is
quasiisomorphic to the cone $\overline{E}^1_\idot$ of the map $V$ of
\eqref{f.v.dia}. This suggests a very useful alternative description
of the category $\DM(G,R)$. Namely, recall that the {\em Tate
  cohomology complex} $\vC^\hdot(G,E_\idot)$ with coefficients in a
complex $E_\idot$ of $R[G]$-modules is by definition the cone of the
trace map \eqref{tr.eq}. Then the map $F$ in \eqref{f.v.dia} induces
a natural map
\begin{equation}\label{phi.ta}
\phi:\overline{E}^1_\idot \to \vC^\hdot(G,E^0_\idot).
\end{equation}
Conversely, given $\overline{E}^1_\idot$ and the map $\phi$, one
recovers $E^1_\idot$ as the cone of the natural map
$$
\begin{CD}
\overline{E}^1_\idot @>{\phi}>> \vC^\hdot(G,E^0_\idot) @>>>
C_\idot(G,E^0_\idot)[1],
\end{CD}
$$
this $E^1_\idot$ comes equipped with the map $V$ of \eqref{f.v.dia},
and $\phi$ then induces the map $F$. We see that $\DM(G,R)$ can be
obtained by inverting quasiisomorphisms in the category of pairs
$\langle E^0_\idot,\overline{E}^1_\idot \rangle$ of complexes of
$R$-modules equipped with an automorphism $\sigma:E^0_\idot \to
E^0_\idot$ of order $p$ and a map $\phi$ of \eqref{phi.ta}.

To rephrase this description even further, let $\wt{P}_\idot$ be the
cone of the augmentation map $P_\idot \to \Z$, where $P_\idot$ is
our fixed projective resolution. Then one can always choose
$E^0_\idot$ to be an $h$-injective complex of $R[G]$-modules, and in
this case, $\vC^\hdot(G,E^0_\idot)$ is quasiisomorphic to the
sum-total complex of the bicomplex
$$
(\wt{P}_\idot \otimes E^0_\idot)^\sigma.
$$
Then the map $\phi$ of \eqref{phi.ta} becomes a $\sigma$-invariant map
$$
\phi:\overline{E}^1_\idot \to \wt{P}_\idot \otimes E^0_\idot,
$$
and the whole data $\langle
E^0_\idot,\overline{E}^1_\idot,\sigma,\phi \rangle$ can be packaged
into a single DG comodule over the DG coalgebra $T_\idot(G,R)$ over
$R$ of the following form:
\begin{equation}\label{t.g}
T_\idot(G,R) = \begin{pmatrix}R[G] & \wt{P}_\idot \otimes R\\0 &
  R\end{pmatrix},
\end{equation}
where $R[G]$ is the group coalgebra of the group $G$, and
$\wt{P}_\idot \otimes R$ is the $R[G]$-comodule obtained from
$\wt{P}_\idot$. The category $\DM(G,R)$ is then equivalent to the
derived category $\D(T_\idot(G,R))$ (that is, the category obtained
by inverting quasiisomorphisms in the category of DG comodules over
$T_\idot(G,R)$, as in \cite[Subsection 1.5.3]{mackey}).

Using DG coalgebras here is essential -- the complex $\wt{P}_\idot$
is infinite, and one cannot simply dualize things and interpret
$\D(T_\idot(G,R))$ as the derived category of DG modules over a DG
algebra. Let us also note that the complex $\wt{P}_\idot \otimes R$
is acyclic, so that $T_\idot(G,R)$ is quasiisomorphic to the
diagonal DG coalgebra $\overline{T}_\idot(G,R)$ with entries $R[G]$
and $R$. However, the derived categories $\D(T_\idot(G,R))$ and
$\D(\overline{T}_\idot(G,R))$ are different: the latter is the
direct sum of $\D(R)$ and $\D(R[G])$, while the former is obtained
by gluing these two categories along the gluing functor given by
Tate cohomology.

\subsection{Maximal Tate cohomology.}\label{tate.subs}

To extend the description of derived Mackey functors in terms of
geometric fixed points and Tate cohomology to arbitrary finite
groups, note that for any group $G$ and derived Mackey functor $E
\in \DM(G,R)$, the value $E([G/\{e\}])$ of $E$ at the free orbit
$[G/\{e\}] \in \Gamma_G$ is acted upon by $G = \Aut_G([G/\{e\}])$
and gives an object in the derived category $\D(R[G])$. For any
subgroup $H \subset G$ with $W_H = N_H/H$, we can apply this
observation to the group $W_H$. Then sending $E \in \DM(G,R)$ to
$\Phi^H(E)([W_H/\{e\}]) \in \D(R[W_H])$ gives a natural functor
$$
\bPhi^H:\D(G,R) \to \D(R[W_H]),
$$
a restricted version of the geometric fixed points functor of
\eqref{infl.h.der}. It has been proved in \cite[Lemma 7.2]{proma}
that this functor has a right-adjoint restricted inflation functor
\begin{equation}\label{binfl}
\bInfl^H:\D(R[W_H]) \to \D(G,R).
\end{equation}
Taking the direct sum of the functors $\bPhi^H$ over all conjugacy
classes of subgroups $H \subset G$, one obtains the functor
\begin{equation}\label{bphi}
\bPhi^\hdot:\DM(G,R) \to \prod_{H \subset G}\D(R[W_H]),
\end{equation}
and one proves the following result.

\begin{lemma}[{{\cite[Lemmas 7.1, 7.3]{proma}}}]\label{ortho.le}
The functor $\bPhi^\hdot$ of \eqref{bphi} is conservative, and for
any $H \subset G$, the restricted inflation functor \eqref{binfl} is
fully faithful. Moreover, let $\DM_H(G,R) \subset \DM(G,R)$ be the
image of the fully faithful embedding $\bInfl^H$. Then every object
$E \in \DM(G,R)$ is an iterated extensions of objects $E_H \in
\DM_H(G,R)$, $H \subset G$, and we have $\Hom(E_H,E_{H'}) = 0$
unless $H$ is conjugate to a subgroup in $H'$.
\end{lemma}

Thus the category $\DM(G,R)$ admits a semiorthogonal decomposition
whose graded pieces $\DM_H(G,R)$ are numbered by conjugacy classes
of subgroup $H \subset G$, and we have $\DM_H(G,R) \cong
\D(R[W_H])$. For any two subgroups $H' \subset H \subset G$, we then
have the gluing functor
\begin{equation}\label{e.h.h}
E^H_{H'} = \bPhi^H \circ \bInfl^{H'}:\D(R[H']) \to \D(R[H]).
\end{equation}
To compute these gluing functors, one introduces in \cite[Subsection
  7.2]{proma} a certain generalization of Tate cohomology that we
call {\em maximal Tate cohomology}. Let us describe it.

Recall that for any subgroup $H \subset G$, the restriction
functor
$$
r^G_H:\D(R[G]) \to \D(R[H])
$$
has a left and right-adjoint induction functor $i^H_G:\D(R[H]) \to
\D(R[G])$.

\begin{defn}[{{\cite[Definition 7.4]{proma}}}]\label{max.tate}
A $\Z[G]$-module $E$ is {\em induced} if $E$ is a direct summand of
a sum of objects of the form $i^H_G(E_H)$, $H \subset G$ a proper
subgroup, $E_H$ a finitely generated $\Z[H]$-module projective over
$\Z$. For any $E \in \D^b(\Z[G])$, the {\em maximal Tate cohomology}
$\vH^\hdot(G,E)$ is given by
$$
\vH^\hdot(G,E) = \RHom^\hdot_{\bD(\Z[G])}(\Z,E),
$$
where $\D_i^b(\Z[G]) \subset \D^b(\Z[G])$ is the smallest
Karoubi-closed triangulated subcategory containing all induced
modules, and $\bD(\Z[G]) = \D^b(\Z[G])/\D^b_i(\Z[G])$ is the
quotient category.
\end{defn}

If $G = \Z/p\Z$, a $\Z[G]$-module is induced if and only if it is
finitely generated and projective, so that maximal Tate cohomology
coincides with the usual Tate cohomology. In general, they are
different. For example, if $G = \Z/n\Z$ is a cyclic group, then
$\vH^\hdot(G,-)=0$ unless $n$ is a prime (this is \cite[Lemma
  7.15~\thetag{ii}]{mackey}).

To compute maximal Tate cohomology, and to generalize it to possibly
infinite coefficients, it is convenient to introduce the following.

\begin{defn}
A complex $P_\idot$ of $\Z[G]$-modules is {\em maximally adapted} if
\begin{enumerate}
\item $P_i=0$ for $i < 0$, $P_0 \cong \Z$, and $P_i$ is induced for
  any $i \geq 1$, and
\item $r^H_G(P_\idot)$ is contractible for any proper subgroup $H
  \subset G$.
\end{enumerate}
\end{defn}

Then for any ring $R$ and maximally adapted complex $P_\idot$, one
defines the Tate cohomology object with coefficients in $E \in
\D(R[G])$ as
\begin{equation}\label{vC}
\vC^\hdot(G,E) = \lim_{\overline{l}{\to}}C^\hdot(G,E \otimes
F^lP_\idot),
\end{equation}
where $F^lP_\idot$ stands for the stupid filtration. It has been
proved in \cite[Subsection 7.2]{proma} that as an object in $\D(R)$,
this does not depend on the choice of $P_\idot$, and for any $E \in
\D^b(\Z[G])$, the complex $\vC^\hdot(G,E)$ computes the maximal Tate
cohomology groups $\vH^\hdot(G,E)$ of Definition~\ref{max.tate}. If
$G$ is a normal subgroup in a larger group $G'$, then for any $E \in
\D(R[G'])$ and an appropriate choice of an adapted complex, the
expression \eqref{vC} also defines $\vC^\hdot(G,E)$ as an object in
$\D(R[G'/G])$.

A full description of the gluing functors \eqref{e.h.h} in terms of
the maximal Tate cohomology objects \eqref{vC} is given in
\cite[Proposition 7.10]{proma}. Here we only reproduce the answer in
the key case $H' = \{e\}$, $H = G$. In this case, we have
$$
E^G_{\{e\}} \cong \vC^\hdot(G,-).
$$
Based on this, in \cite[Section 6]{mackey}, one develops a description of
the category $\DM(G,R)$ in terms of a certain upper-triangular DG
coalgebra $T_\idot(G,R)$ similar to the coalgebra
\eqref{t.g}. Unfortunately, this is rather heavy technically (in
particular, one cannot get a genuine DG coalgebra and has to settle
for its $A_\infty$-version).

A big simplification occurs in the case of the cyclic group $G =
\Z/n\Z$, $n \geq 1$. As we have mentioned, $vH^\hdot(\Z/n\Z,-)=0$
unless $n$ is a prime; moreover, for any two primes $p$, $p'$, it
has been proved in \cite[Lemma 9.10]{proma} that
$$
\vC^\hdot(\Z/p\Z,\vC(\Z/p'\Z,-)) = 0,
$$
so that the compositions of the gluing functors \eqref{e.h.h} vanish
tautologically. Because of this, it has been possible to obtain a
reasonably simple description of $\DM(G,R)$. Let us reproduce it
(this is \cite[Subsection 9.4]{proma}).

Fix an integer $n \geq 1$, let $G = \Z/n\Z$, and denote by $I_n$ the
groupoid of $G$-orbits. These are numbered by divisors of $n$, so
that explicitly, we have
$$
I_n = \coprod_{d | n} \ppt_d,
$$
where $\ppt_d$ stands for the groupoid with one object with
automorphism group $\Z/d\Z$. For any prime $p$, let
$$
I^p_n = \coprod_{p | d | n} \ppt_d \subset I_n,
$$
and let $I^\hdot_n$ be the disjoint union of $I^p_n$ over all primes
$p$. Of course $I^p_n$ is empty unless $p$ divides $n$. For any $p$
dividing $n$, we have the natural embedding $i:I^p_n \to I_n$, and
we also have a natural projection $\pi:I^p_n \to I_n$ given by the
union of quotient maps $\ppt_{pm} \to \ppt_m$, $m$ a divisor of
$n/p$. Let
\begin{equation}\label{i.pi}
i,\pi:I^\hdot_n \to I_n
\end{equation}
be the disjoint union of these functors. Finally, choose a
projective resolution $P_\idot$ of the constant functor $\Z \in
\Fun(I^\hdot_n,\Z)$, and let $\wt{P}_\idot$ be the cone of the
augmentation map $P_\idot \to \Z$.

\begin{defn}
An {\em $R$-valued fixed points datum} for $I_n$ is a pair $\langle
M_\idot,\alpha \rangle$ of a complex $M_\idot$ in the category
$\Fun(I_n,R)$ and a map
$$
\alpha:\pi^*M_\idot \to \wt{P}_\idot \otimes i^*M_\idot,
$$
where $i$ and $\pi$ are the projections \eqref{i.pi}.
\end{defn}

Then $I_n$-fixed points data form a category, and inverting
quasiisomorphisms in this category, one obtains a category
$\DA_n(R)$. It has been shown in \cite[Subsection 9.3]{proma} that
$\DA_n(R)$ is a triangulated category that does not depend on the
coice of a resolution $P_\idot$. One then proves the following
comparison result.

\begin{prop}[{{\cite[Proposition 9.14~\thetag{i}]{proma}}}]\label{Z.n.prop}
For any $n \geq 1$ and ring $R$, there exiss a natural equivalence
of triangulated categories
$$
\DM(\Z/n\Z,R) \cong \DA_n(R).
$$
\end{prop}

\section{Mackey profunctors.}

\subsection{Definitions.}

Assume now given an infinite group $G$. Then if one considers
arbitrary $G$-sets, all the theory of Section~\ref{mack.sec} and
Section~\ref{dm.sec} becomes trivial (for example, the Burnside ring
of \eqref{A.H.eq} would be identically zero). If we stick to finite
$G$-sets, the theory works, up to a point. However, finer parts of
the theory such as Subsection~\ref{tate.subs} break down.

As it turns out, there is an interesting possibility in between the
two extremes. It is based on the following definition.

\begin{defn}[{{\cite[Definition 3.1]{proma}}}]
A $G$-set $S$ is {\em admissible} if
\begin{enumerate}
\item for any $s \in S$, its stabilizer $G_s \subset G$ is a
  cofinite subgroup, and
\item for any cofinite subgroup $H \subset G$, the fixed points set
  $S^H$ is finite.
\end{enumerate}
\end{defn}

Explicitly, any $G$-set $S$ decomposes into a disjoint union of
orbits
\begin{equation}\label{copr.eq}
S = \coprod_{s \in S/G}[G/H_s].
\end{equation}
Then $S$ is admissible iff all the subgroups $H_s \subset G$ are
cofinite, and for any cofinite subgroup $H \subset G$, at most a
finite number of them contain $H$.

One denotes the category of admissible $G$-sets by $\wGamma_G$, and
one notes that this category has fibered products. Therefore we can
consider the category $Q\wGamma_G$ and the embedding \eqref{e.eq}.

\begin{defn}\label{pro.def}
For any ring $R$, an object $E \in \D(\wGamma_G^o,R)$ is {\em
  additive} if for any admissible $G$-set $S$ with decomposition
\eqref{copr.eq}, the natural map
$$
E(S) \to \prod_{s \in S/G}E([G/H_s])
$$
is an isomorphism. An {\em $R$-valued $G$-Mackey profunctor} is an
object $E \in \Fun(Q\wGamma_G,R)$ whose restriction $e^*E$ is
additive.
\end{defn}

We denote by $\wM(G,R) \subset \Fun(Q\Gamma_G,R)$ the full
subcategory spanned by Mackey profunctors. Since infinite sums in
the category of $R$-modules are exact, $\wM(G,R)$ is an abelian
category.

As in the finite group case, for any cofinite subgroup $H \subset
G$, we have the pair of adjoint functors $\rho^H$, $\gamma^H$
between $\wGamma_G$ and $\wGamma_H$. The functors preserve fibered
products, the corresponding functors $Q(\rho^H)$, $Q(\gamma^H)$ are
also adjoint, and the functor $Q(\rho^H)_! \cong Q(\gamma^H)^*$
preserves the additivity condition of Definition~\ref{pro.def}, thus
sends Mackey profunctors to Mackey profunctors. The corresponding
functor
$$
\Psi^H = Q(\rho^H)_! \cong Q(\gamma^H)^*:\wM(G,R) \to \wM(H,R)
$$
is the {\em categorical fixed points} functor. Moreover, if we
denote $W=N_H/H$, then the fixed points functor $\phi^H:\wGamma_G
\to \Gamma_W$ also preserves fibered products, and we have an
adjoint pair of functors
\begin{equation}\label{infl.pro}
\begin{aligned}
\Phi^H &= Q(\phi^H)_!:\Fun(Q\wGamma_G,R) \to \Fun(Q\Gamma_W,R),\\
\Infl^H &= Q(\phi^H)^*:\Fun(Q\Gamma_W,R) \to \Fun(Q\wGamma_G,R).
\end{aligned}
\end{equation}
For the same reasons as in the usual case, these preserve
additivity, thus induce functors between $\wM(G,R)$ and $\M(W,R)$,
the {\em geometric fixed points functor} and the {\em inflation
  functor}.

Lemma~\ref{infl.le} also holds for Mackey profunctors. In
particular, for any normal cofinite subgroup $N \subset G$, the
inflation functor $\Infl^N$ is fully faithful. Any Mackey profunctor
$E$ gives rise to a Mackey functor $E_N = \Phi^NE \in \M(W,R)$, and
we have a natural surjective map
$$
E \to \Infl^N(E_N).
$$
For any two cofinite normal subgroups $N \subset N' \subset G$, we
have a natural isomorphism
\begin{equation}\label{e.n}
\Phi^{N'/N}E_N \cong E_{N'}.
\end{equation}
It is convenient to axiomatize the situation as follows.

\begin{defn}
An {\em $R$-valued $G$-normal system} $\langle E_\idot \rangle$ is a
collection of Mac\-key functors $E_N \in \M(G/N,R)$, one for each
cofinite normal subgroup $N \subset G$, and a collection os
isomorphisms \eqref{e.n}, one for each pair of cofinite normal
subgroups $N \subset N' \subset G$.
\end{defn}

Normal systems form an additive $R$-linear category denoted
$\N(G,R)$. Sending $E$ to $\langle \Phi^N(E) \rangle$ gives a
functor $\Phi:\wM(G,R) \to \N(G,R)$.  This functor has a
right-adjoint
\begin{equation}\label{infl.tot}
\Infl:\N(G,R) \to \wM(G,R)
\end{equation}
sending a normal system $\langle E_N \rangle$ to
\begin{equation}\label{infl.tot.pro}
E = \lim_{\overset{N}{\gets}}\Infl^N(E_N),
\end{equation}
where the limit is taken over all cofinite normal subgroups $N
\subset G$, with respect to the natural maps adjoint to the
isomorphisms \eqref{e.n}. It turns out that the following is true.

\begin{lemma}[{{\cite[Proposition 3.5]{proma}}}]\label{ns.le}
The functor $\Infl$ of \eqref{infl.tot} is fully
faithful. Conversely, for any Mackey profunctor $E \in \wM(G,R)$,
the adjunction map
\begin{equation}\label{can.filt}
E \to \Infl(\Phi(E)) = \lim_{\overset{N}{\gets}}\Infl^N(\Phi^N(E))
\end{equation}
is surjective.
\end{lemma}

One says that a Mackey profunctor $E$ is {\em separated} if the
surjective map \eqref{can.filt} is bijective -- or equivalently, if
$E$ lies in the image of the fully faithful embedding $\Infl$. One
denotes by $\wM_s(G,R) \subset \wM(G,R)$ the full subcategory
spanned by separated Mackey profunctors. This category is additive
but not necessarily abelian. However, it is equivalent to $\N(G,R)$,
and this allows to transport results about Mackey functors to
separated Mackey profunctors. In particular, it has been proved in
\cite[Lemma 3.9]{proma} that the natural embedding $\wM_s(G,R)
\subset \Fun(Q\wGamma_G,R)$ admits a left-adjoint additivization
functor
\begin{equation}\label{add.pro}
\Add:\Fun(Q\wGamma_G,R) \to \wM_s(G,R).
\end{equation}
Using this functor, one can extend the definition of the geomeric
fixed points functor $\Phi^H$ to an arbitrary subrgroup $H \subset
G$. Indeed, even if $H \subset G$ is not cofinite, the inflation
functor $\Infl^H$ of \eqref{infl.pro} sends Mackey profunctors to
Mackey profunctors, and it also sends separated Mackey profunctors
to separated ones. Then it has a natural left-adjoint given by
$$
\Phi^H = \Add \circ Q(\phi^H)_!.
$$
Another application of the addivization functor is tensor products
of separated Mackey profunctors; these are defined by the same
formula \eqref{ma.produ} as in the Mackey functor case. The product
is associative, commutative and unital. The unit object is the
Burnside Mackey profunctor $\wA$. Explicitly, it is given by
$$
\wA([G/H]) = \wh{A}^H,
$$
where the completed Burnside ring $\wh{A}^H$ is given by
\eqref{A.H.eq} with $\Gamma_H$ replaced by $\wGamma_H$. By virtue of
the decomposition \eqref{copr.eq}, we have
\begin{equation}\label{a.g.pro}
\wh{A}^G = \prod_{H \subset G}\Z,
\end{equation}
where the product is over conjugacy classes of cofinite subgroups $H
\subset G$.

\subsection{An example.}

To illustate the difference between Mackey functors and Mackey
profunctors, consider the case $G=\Z$, the infinite cyclic
group. Cofinite subgroups in $\Z$ are of the form $n\Z \subset \Z$,
$n \geq 1$, so that $\Z$-orbits are numbered by positive
integers. The Burnside ring $A^\Z$ is given by \eqref{a.g};
explicitly, as an abelian group, we have
$$
A^\Z = \Z[\eps_1,\eps_2,\dots],
$$
where the generators $\eps_n$, $n \geq 1$ correspond to $\Z$-orbits
$[\Z/n\Z]$. For any $n,m \geq 1$, the product $[\Z/m\Z] \times
[\Z/n\Z]$ is the union of copies of the orbit $\Z/\{n,m\}\Z$, where
$\{n,m\}$ is the least common multiple of $n$ and $m$. Since the
cardinality of this product is $nm$, this implies that
\begin{equation}\label{eps.n.m}
\eps_n \eps_m = \frac{nm}{\{n,m\}}\eps_{\{n,m\}},
\end{equation}
and this completely defines the product in the Burnside ring $A^\Z$.

Now, for the completed Burnside ring $\wh{A}^\Z$, the sum
\eqref{a.g} is replaced by the product \eqref{a.g.pro}. Therefore we
have
$$
\wh{A}^\Z = \Z\{\eps_1,\eps_2,\dots\},
$$
that is, the group of infinite linear combinations of the generators
$\eps_n$, $n \geq 1$. The product is still given by \eqref{eps.n.m}.

Here is one observation we can make right away: the completed
Burnside ring $\wh{A}^\Z$ is isomorphic to the universal Witt
vectors ring $\W(\Z)$. There is a reason for this coincidence, but
it goes beyond the subject of the present paper.

Another observation is the following. For any ring $R$, any
$R$-valued $\Z$-Mackey functor $E$ is acted upon by the Burnside
ring $A^\Z$, and any $\Z$-Mackey profunctor $E$ is acted upon by the
completed Burnside ring $\wh{A}^\Z$. Assume that $R$ is $p$-local
for some prime $p$ --- that is, every integer $n$ prime to $p$ is
invertible in $R$. Then for any such $n$, we have a well-defined
endomorphism
\begin{equation}\label{eps.n}
\frac{1}{n}\eps_n:E \to E,
\end{equation}
and \eqref{eps.n.m} shows that this endomorphism is idempotent. Thus
$E$ is naturally equipped with a large family of commuting
idempotents.

If $E$ is a $\Z$-Mackey functor, then this is the end of the
story. However, if $E$ is a $\Z$-Mackey profunctor, we can replace
the commuting idempotents \eqref{eps.n} with orthogonal commuting
idempotents $\eps_{(n)}$ given by
$$
\eps_{(n)} = \frac{1}{n}\eps_n \cdot \prod_{i
  \text{ does not divide }n}\left(1 - \frac{1}{i}\eps_i\right),
$$
where the product is over $i$ prime to $p$. We have
$$
1 = \sum_{n \text{ prime to }p}\eps_{(n)},
$$
so that for any $E \in \wM(\Z,R)$, we have a canonical decomposition
\begin{equation}\label{typ.deco}
E = \prod_{n \text{ prime to } p}E_{(n)}, \qquad E_{(n)} = \Im
\eps_{(n)} \subset E.
\end{equation}
This is functorial in $E$ and gives a decomposition of the category
$\wM(\Z,R)$.

It turns out that the pieces of this decomposition can be described
in terms of Mackey profunctors for the group $\Z_p$ of $p$-adic
integers. Namely, by definition, the category $\wM(G,R)$ only
depends on the profinite completion of the group $G$. The profinite
completion $\wh{\Z}$ is given by
$$
\wh{\Z} = \prod_l \Z_l,
$$
where the product is over all primes $l$. If we denote by $\Z'_p$
the product of all primes different from $p$, then $\wh{\Z} = \Z_p
\times \Z'_p$, and we have the geomeric fixed points functor
$$
\Phi^{(p)} = \Phi^{\Z'_p}:\wM(\Z,R) \cong \wM(\wh{\Z},R) \to
\wM(\Z_p,R).
$$
In \cite[Subsection 9.2]{proma}, this functor has been refined --
for any integer $n \geq 1$ prime to $p$, one constructs a functor
$$
\Phi^{(p)}_{[n]}:\wM(\Z,R) \to \wM(\Z_p,R[\Z/n\Z])
$$
such that $\Phi^{(p)}=\Phi^{(p)}_{(1)}$, and one proves the
following result.

\begin{prop}[{{\cite[Proposition 9.4]{proma}}}]
Assume that the ring $R$ is $p$-local. Then the functor
$$
\prod \Phi^{(p)}_{(n)}:\wM(\Z,R) \to \prod_{n \text{ prime to
  }p}\wM(\Z_p,R[\Z/n\Z])
$$
is an equivalence of categories, and for any $E \in \wM(\Z,R)$, the
component $E_{(n)}$ of the decomposition \eqref{typ.deco}
corresponds to $\Phi^{(p)}_{(n)}(E)$.
\end{prop}

If one identifies the Burnside ring $\wh{A}^\Z$ with the Witt
vectors ring $\W(\Z)$, then the decomposition \eqref{typ.deco}
corresponds to the $p$-typical decomposition of Witt vectors;
because of this, in \cite[Subsection 9.2]{proma}, \eqref{typ.deco}
is called {\em $p$-typical decomposition}.

\subsection{Derived version.}

The derived counterpart of the theory of Mackey profunctors is
largely parallel to the theory of derived Mackey functors of
Section~\ref{dm.sec}. Since the category $\wGamma_G$ has fibered
products, it can be plugged into the machinery of
Subsection~\ref{sc.subs}. We thus have the category $S\wGamma_G$ and
the category $\DS(\wGamma_G,R)$ for any ring $R$. The additivity
condition of Definition~\ref{pro.def} makes sense of objects in the
derived category $\D(\wGamma_G^o,R)$. A {\em derived $R$-valued
  $G$-Mackey profunctor} is an object $E \in \DS(\wGamma_G,R)$ whose
base part $\overline{E}$ of \eqref{bp.eq} is additive.

The category of derived $G$-Mackey profunctors is denoted by
$\wDM(G,R)$. The standard $t$-structure on $\DS(\wGamma_G,R)$
induces a $t$-structure on $\wDM(G,R)$, and the heart of this
$t$-structure is identified with $\wM(G,R)$.

For any cofinite subgroup $H \subset G$, we have the categorical
fixed points functor
$$
\Psi^H:\wDM(G,R) \to \wDM(H,R)
$$
defined exactly as in \eqref{psi.h.der}. Moreover, if we denote $W =
N_H/H$, then we have the adjoint pair
$$
\Phi^H:\wDM(G,R) \to \DM(W,R), \quad \Infl^H:\DM(W,R) \to \wDM(G,R)
$$
of the geometric fixed points functor $\Phi^H$ and the inflation
functor $\Infl^H$ defined as in \eqref{infl.h.der}. The inflation
functor is defined even for a subgroup $H \subset G$ that is not
cofinite. For a cofinite normal subgroup $N \subset G$, the
inflation functor $\Infl^N$ is fully faithful (this is \cite[Lemma
  6.12]{proma}). We also have a version of Lemma~\ref{infl.der.le}
(the statement is the same, but one needs to use a slightly stronger
notion of an adapted set given in \cite[Definition 6.3]{proma}).

To proceed further, one needs to introduce a derived counterpart of
the notion of a normal system. Simply repeating
Definition~\ref{ns.def} does not work since it does not produce a
triangulated category -- we have to package the same data in a more
elaborate way. To do this, let $\Nn$ be the partially ordered set of
normal cofinite subgroups $N \subset G$, ordered by inclusion, and
let
$$
\hGamma_G \subset \Gamma_G \times \Nn^o
$$
be the subcategory of pairs $\langle S,N \rangle$ such that $N$ acts
trivially on $S$ (we treat the partially ordered set $\Nn$ as a
small category in the usual way, and let $\Nn^o$ be the opposite
category). For every $N \in \Nn$, we have a natural functor
$\tau_N:\Gamma_{G/N} \to \hGamma_G$ sending $S$ to $\langle S,N
\rangle$. Moreover, for any pair of cofinite normal subgroups $N
\subset N' \subset G$, we have the fixed points functor
$\phi^{N'/N}:\Gamma_{G/N} \to \Gamma_{G/N'}$, and the inclusions
$S^{N'/N} \subset S$ glue together to give a map of functors
\begin{equation}\label{tau.n.n}
\tau_{N'} \circ \phi^{N'/N} \to \tau_N.
\end{equation}
We also have the forgetful functor $\nu:\hGamma_G \to \Nn^o$ sending
$\langle S,N \rangle$ to $N \in \Nn^o$, and to define derived normal
systems, we do the $S$-construction of Subsection~\ref{sc.subs}
relatively over $\Nn^o$ -- that is, we consider the full subcategory
$$
S(\hGamma_G/\Nn^o) \subset S\hGamma_G
$$
formed by diagrams \eqref{n-ex} in $\hGamma_G$ such that all maps
become invertible after applying $\nu$. Then for any $N \in \Nn^o$,
the functor $\tau_N$ gives a functor
$$
S(\tau_N):S\Gamma_{G/N} \to S(\hGamma_G/\Nn^o),
$$
and for a pair $N \subset N' \subset G$, $N,N' \in \Nn$, the
morphism \eqref{tau.n.n} gives a morphism
\begin{equation}\label{S.tau.n.n}
S(\tau_N) \to S(\tau_{N'}) \circ S(\phi^{N'/N}).
\end{equation}
Therefore any object $E \in \D(S(\hGamma_G/\Nn),R)$ produces a
collection of objects
$$
E_N = S(\tau_N)^*E \in \D(S\Gamma_{G/N},R), \qquad N \in \Nn,
$$
related by morphisms
\begin{equation}\label{e.n.bis}
E_N \to S(\phi^{N'/N})^*E_{N'}.
\end{equation}

\begin{defn}\label{ns.def}
An {\em $R$-valued derived $G$-normal system} is an object $E \in
\D(S(\hGamma_G/\Nn),R)$ such that for any $N \in \Nn$, $E_N \in
\D(S\Gamma_{G/N},R)$ is a Mackey functor, and for any $N \subset N'
\subset G$, the map $\Phi^{N'/N}E_N \to E_{N'}$ adjoint to the map
\eqref{e.n.bis} is an isomorphism.
\end{defn}

By construction, derived normal systems form a triangulated
category; we denote it by $\DN(G,R)$. For every integer $n$, we let
$\DN^{\leq n}(G,R) \subset \DN(G,R)$ be the full subcategory formed
by objects $E$ such that for any $N \in \Nn$, $E_N$ lies in
$\DM^{\leq n}(G/N,R)$. These subcategories do not necessarily give a
$t$-structure on $\DN(G,R)$. However, they are perfectly
well-defined.

With these definitions, one constructs a functor
$$
\Phi:\wDM(G,R) \to \DN(G,R)
$$
such that for any $E \in \DM(G,R)$ and $N \in \Nn$, $\Phi(E)_N$ is
canonically identified with $\Phi^N(E)$, and the map \eqref{e.n} is
adjoint to the map \eqref{e.n.bis}. Then one proves the following
derived counterpart of Lemma~\ref{ns.le}.

\begin{prop}[{{\cite[Proposition 8.2]{proma}}}]\label{ns.prop}
Assume that the group $G$ is finitely generated. Then for every
integer $n$, the functor $\Phi$ induces an equivalence of categories
$$
\Phi:\DM^{\leq n}(G,R) \cong \DN^{\leq n}(G,R).
$$
\end{prop}

As a corollary of this Proposition, we can consider the union
$$
\DN^-(G,R) = \bigcup_n \DN^{\leq n}(G,R) \subset \DN(G,R)
$$
of all the categories $\DN^{\leq n}(G,R)$, and see that it is
naturally equivalent to the full subcategory $\wDM^-(G,R) \subset
\wDM(G,R)$ of derived Mackey profunctors bounded from above with
respect to the standard $t$-structure.

As we see, Proposition~\ref{ns.prop} is much stronger that
Lemma~\ref{ns.le}: a derived Mackey profunctor is separated as soon
as it is bounded from above. In particular, any Mackey profunctor $E
\in \wM(G,R) \subset \wDM(G,R)$ is separated in the derived
sense. What happens is, the equivalence $\Infl$ inverse to $\Phi$ is
given explicitly by
$$
\Infl(E) = \dlim_{\overset{N}{\gets}}\Infl^N(E_N),
$$
but since the inverse limit is not an exact functor, this is
different from \eqref{infl.tot.pro} -- there could be non-trivial
contributions from $R^1\lim_{\gets}$ in the right-hand side. One
shows that we in fact have $R^i\lim_\gets = 0$ for $i \geq 2$, so
for any Mackey profunctor $E$, the kernel of the canonical
surjective map \eqref{can.filt} is identified with
$$
R^1\lim_{\overset{N}{\gets}}\Infl^N(L^1\Phi^N(E)).
$$
Applications of Proposition~\ref{ns.prop} are similar to those of
Lemma~\ref{ns.le}. Firstly, one proves in \cite[Lemma 8.3]{proma}
that the embedding
$$
\wDM^-(G,R) \subset \DS^-(\wGamma_G,R)
$$
admits a left-adjoint additivization functor
\begin{equation}\label{add.pro.der}
\Add:\DS^-(\wGamma_G,R) \to \wDM^-(G,R).
\end{equation}
We note that by adjunction, $\Add$ is right-exact with respect to
the standard $t$-structures, and on the hearts of the standard
$t$-structures, it induces a functor $\Fun(Q\wGamma_G,R) \to
\wM(G,R)$ right-adjoint to the natural embedding --- that is, a
refinement of the additivization functor \eqref{add.pro}. Using this
refinement, one extends the product \eqref{ma.produ} to all Mackey
profunctors, not only the separated ones, and one does the same with
the geometric fixed points functors $\Phi^H$ for arbitrary subgroups
$H \subset G$.

In the derived theory, one uses \eqref{add.pro.der} to define the
tensor product of derived Mackey profunctors by
\eqref{ma.produ.der}, and one defines the geometric fixed points
functor $\Phi^H$ with respect to an arbitrary subgroup $H \subset G$
by
$$
\Phi^H = \Add \circ S(\phi^H)_!.
$$
This is left-adjoint to the inflation functor $\Infl^H$. We also
have a version of Lemma~\ref{sph.le} --- namely, \cite[Lemma
  8.7]{proma} --- and can prove an analog of Lemma~\ref{ortho.le},
although the semiorthognal decomposition on $\DM^-(G,R)$ would have
an infinite number of terms.

Finally, let us note that in the case $G=\Z$, the infinite cyclic
group, we have a complete analog of
Proposition~\ref{Z.n.prop}. Namely, let
$$
I = \coprod_{n \geq 1}\ppt_n
$$
be the groupoid of $\Z$-orbits, let $I^p \subset I$ be the
subcategory spanned by $\ppt_{np}$, $n \geq 1$, let $I^\hdot$ be the
disjoint union of the categories $I^p$ over all prime $p$, and let
\begin{equation}\label{i.pi.bis}
i,\pi:I^\hdot \to I
\end{equation}
be the natural functors defined in the same way as
\eqref{i.pi}. Choose a projective resolution $P_\idot$ of the
constant functor $\Z \in \Fun(I^\hdot,\Z)$, and let $\wt{P}_\idot$
be the cone of the augmentation map $P_\idot \to \Z$.

\begin{defn}
An {\em $R$-valued fixed points datum} is a pair $\langle
M_\idot,\alpha \rangle$ of a complex $M_\idot$ in the category
$\Fun(I,R)$ bounded from above, and a map
$$
\alpha:\pi^*M_\idot \to \wt{P}_\idot \otimes i^*M_\idot,
$$
where $i$ and $\pi$ are the projections \eqref{i.pi.bis}.
\end{defn}

Then just as in the case of a finite cyclic group, $I$-fixed points
data form a category, and inverting quasiisomorphisms in this
category, one obtains a category $\DA(R)$. One shows that $\DA(R)$
is a triangulated category that does not depend on the choice of a
resolution $P_\idot$, and proves the following result.

\begin{prop}[{{\cite[Proposition 9.14~\thetag{ii}]{proma}}}]
For any ring $R$, there exiss a natural equivalence of triangulated
categories
$$
\DM^-(\Z,R) \cong \DA(R).
$$
\end{prop}

\bigskip

\noindent
{\sc
Steklov Math Institute, Algebraic Geometry section\\
Laboratory of Algebraic Geometry, NRU HSE\\
\mbox{}\hspace{30mm}and\\
IBS Center for Geometry and Physics, Pohang, Rep. of Korea}

\bigskip

\noindent
{\em E-mail address\/}: {\tt kaledin@mi.ras.ru}

\end{document}